\documentclass[a4paper,11pt]{article}

\usepackage{graphicx,amssymb,amsthm}
\usepackage{epsf}

\columnwidth\textwidth

\def\al{\alpha}
\def\be{\beta}
\def\ga{\gamma}

\def\si{\sigma}

\def\Ga{\Gamma}

\def\map{\rightarrow}
\def\bq{\begin{equation}}
\def\eq{\end{equation}}

\def\rc{{\mathbb R}}

\def\ti{\tilde}

\def\Ga{\Gamma}

\def\ss{\subset}

\def\ot{\otimes}
\def\zv{\bfseries\itshape}
\def\bw{\bigwedge}
\def\ca{{\cal A}}
\def\dvapiz{2\pi\sqrt{-1}{\mathbb Z}}

\begin{document}

\begin{center}
{\Large \scshape Quantization of Poisson families and of twisted
Poisson structures}

\vskip 7mm

Pavol \v Severa\\{\it Dept.~of Theoretical Physics\\Comenius
University\\Bratislava, Slovakia}\vskip 7mm
\end{center}

\section*{Introduction}
 If $\phi$ is a closed 3-form on a smooth manifold $M$, a $\phi$-Poisson structure is, according to \cite{sw}, a
 bivector field $\pi$ satisfying the equation
 $$[\pi,\pi]=2\wedge^3\tilde\pi(\phi),$$
where $\ti\pi:T^*M\map TM$ is the linear map given by $\pi$. More
conceptually, it is a Dirac structure in an exact Courant
algebroid (see the appendix for the definitions). In this paper we
try to quantize this structure if $\pi$ depends formally on an
indeterminate $\hbar$ and $\pi=O(\hbar)$; $\phi$ is also allowed
to depend on $\hbar$. Under the condition that the periods of
$\phi$ are in $\dvapiz$ (here $\pi$ is $3.14\dots$, not a
bivector) we construct a stack of algebras on $M$; the stack
depends on a choice of an element $\Phi\in H^3(M,\dvapiz)$  and of
a Dirac structure in an exact Courant algebriod $E$, such that the
image of $\Phi$ in $H^3(M,{\mathbb C})$ is the characteristic
class of $E$. The stack can reasonably be called a deformation of
the ``gerbe" corresponding to $\Phi$.

The basic idea is as follows: since locally the cohomology class
of $\phi$ vanishes, $\pi$ is locally equivalent to a Poisson
structure. This local Poisson structure is not unique, but any two
are connected by a diffeomorphism formally depending on $\hbar$,
and hence they give rise to isomorphic $*$-product. To keep track
of these diffeomorphisms and isomorphisms we introduce a kind of
families of Poisson structures, called tight families here; they
are given by a Maurer-Cartan equation, and hence their
quantization follows immediately from Kontsevich's Formality
theorem \cite{def}. The idea and the interpretation of these
families is again due to Kontsevich \cite{ko}.

\subsection*{Acknowlegdement}
I am particularly grateful to P.~Bressler and P.~Xu who both
suggested me that Dirac structures in exact Courant algebroids
should be quantized to stacks, and to B.~Jur\v co for showing me
his independent and very similar results. I would also like to
thank to J.~Stasheff for interesting comments.

\section{Tight families of associative algebras}

This section follows closely the appendix A.2 of \cite{ko}, it is
here just for reader's convenience.  Let $\ca$ be a vector space
with a marked vector $1_{\ca}\in \ca$, $B$ a manifold and $\chi$ a
closed 3-form on $B$; a {\zv $\chi$-tight family of algebras
indexed by $B$} is an element $\ga\in\Ga(\bw T^*B)\otimes
C^\bullet(\ca,\ca)$ of total degree 2 (where $C^k(\ca,\ca)$ is the
space of $k$-linear maps from $\ca$ to $\ca$), satisfying a
modified Maurer--Cartan equation\footnote{if we augment
$C^\bullet(\ca,\ca)$ by setting $C^{-1}(\ca,\ca)=\mathbb R$ (or
$\mathbb C$) and $d\,1=1_\ca$, this becomes the ordinary
Maurer-Cartan equation}
$$d\ga+{1\over 2}[\ga,\ga]=\chi\otimes 1_\ca$$
and some additional constraints involving $1_\ca$ (see below). In
terms of bihomogeneous components ($\ga=\ga_0+\ga_1+\ga_2$, the
bidegree of $\ga_k$ is $(k,2-k)$)

\begin{enumerate}
\item $[\ga_0,\ga_0]=0$, i.e.~$\ga_0$ is a family of associative
algebras (on the fixed vector space $\ca$) indexed by $B$,
\item $[\ga_1,\ga_0]+d\ga_0=0$, i.e.~if $\ga_1$ is interpreted as
a connection on the vector bundle $\ca\times B\map B$, the algebra
structure  on the fibres is parallel,
\item $[\ga_2,\ga_0]+([\ga_1,\ga_1]/2+d\ga_1)=0$, i.e.~the
connection $\ga_1$ need not be flat, but the result of a parallel
transport along an infinitesimal closed curve is an inner
automorphism given by $\ga_2$,
\item $[\ga_2,\ga_1]+d\ga_2=\chi\otimes 1_\ca$;
\end{enumerate}

\noindent the additional constraints are
\begin{enumerate}
\setcounter{enumi}{4}
\item $1_\ca$ is the unit for every algebra,
\item $1_\ca$ is parallel.
\end{enumerate}

Suppose that $\ca$ is finite-dimensional, so that parallel
transport makes sense. Then for any curve in $B$ the parallel
transport gives us an isomorphism between the algebras at the
endpoints; moreover, for a closed curve spanned by a disk $D$ this
automorphism is inner, given by an element $a_{D,x}\in\ca$ (to be
constructed below), where $x$ is the starting point of the curve.
If moreover $D'$ is another disk homotopic to $D$ rel boundary
then $a_{D',x}=a_{D,x}\exp\int\chi$, where the integral is over
the volume swept by the homotopy (this is the meaning of 4.).
These $a$'s behave in the obvious nice way under composition of
curves and disks:
$$\epsfxsize=7cm \epsfbox{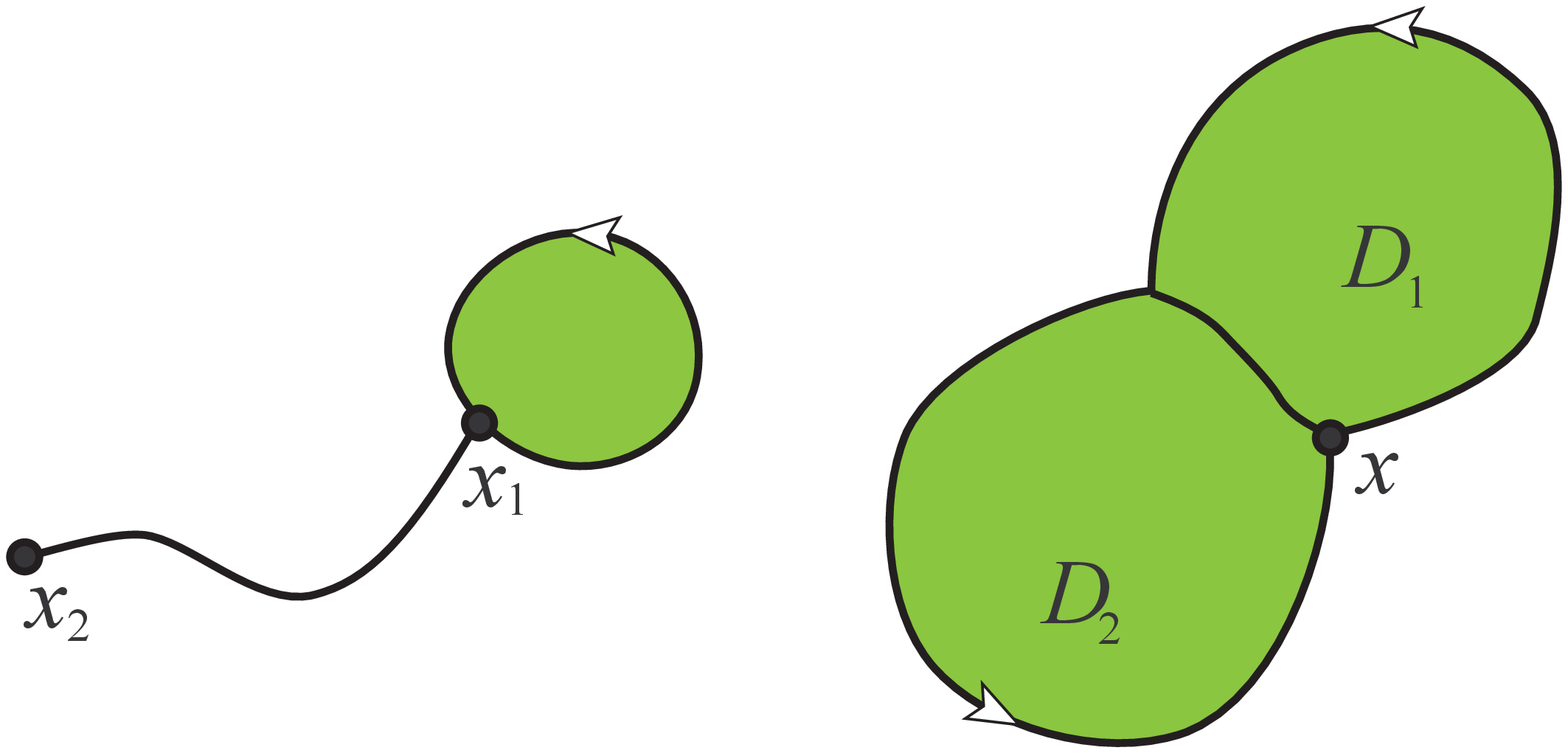}$$

\noindent For the left picture, if we add to our closed curve a
path from $x_1$ to $x_2$ and back again (with no other changes to
the disk) then $a_{D,x}$ doesn't change for any $x$ lying on the
original closed curve; on the other hand,
$a_{D,x_2}=T(a_{D,x_1})$, where $T$ is the parallel transport from
$x_1$ to $x_2$. For the right picture, if the disk is decomposed
to $D_1$ and $D_2$ then $a_{D,x}=a_{D_1,x}a_{D_2,x}$.

These two composition properties also make it clear how to
actually construct $a_{D,x}$:
$$\epsfxsize=4.5cm \epsfbox{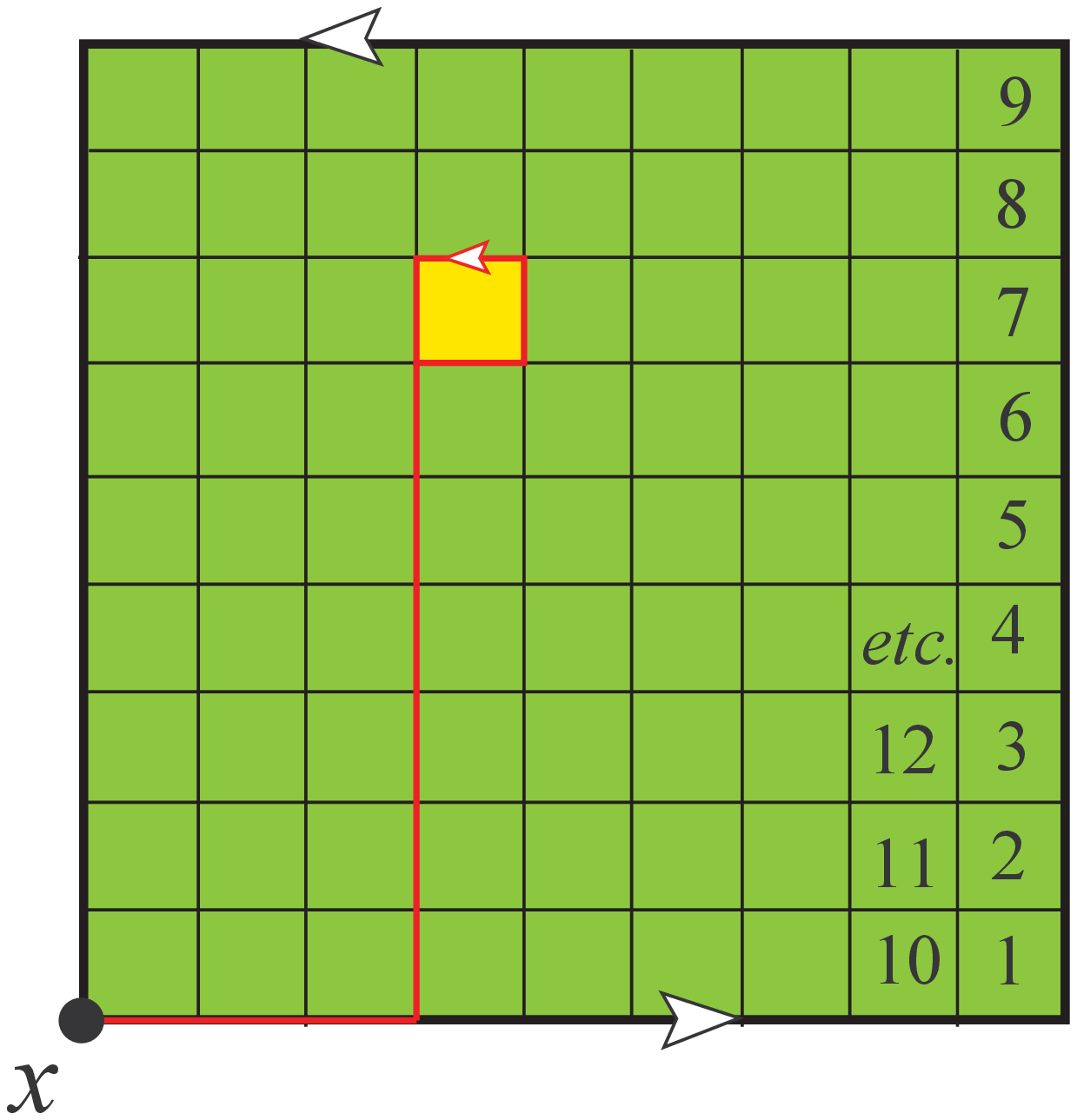}$$
For example, decompose $D$ as the square on the picture and number
the little squares as indicated, then connect each of them with
$x$ by a line that is first horizontal and then vertical. In this
way we express $a_{D,x}$ as a product of $a$'s for the little
squares. Hence it is enough to know $a_D$'s for infinitesimal
$D$'s and these are (by 3.) $\exp(\int_D\ga_2)$.

Finally, there are two natural group actions on tight families. If
$\be$ is a 2-form on $B$, we can add $\be\otimes 1_\ca$ to $\ga$,
getting a $\chi+d\be$-tight family; we shall call it {\zv outer
transformation}. Notice that only $\ga_2$ gets changed; the
connection remains the same and $a_{D,x}$'s are multiplied by
$\exp(\int_D\be)$. Secondly, the Lie algebra of elements
$\al\in\Ga(\bw T^*B)\otimes C^\bullet(\ca,\ca)$ of total degree 1
acts on the space of $\chi$-tight families by
$\al\cdot\ga=d\al+[\al,\ga]$, we shall call this {\zv inner
transformations}. In particular, the outer transformations by
exact 2-forms are also inner.

(Let us finish with an unimportant remark that in a ``global
version'' of tight families one should have a line bundle over the
loop space of $B$, coming from a gerbe on $B$; $a_{D,x}$'s  should
be replaced by $a_p$'s, where $p$'s are elements of the associated
principal bundle.)

\section{Tight families of $*$-products}

A {\zv $\chi$-tight family of $*$-products on $M$ indexed by $B$}
appears when we take $\ca=C^\infty(M)[[\hbar]]$ and instead of
$C^\bullet(\ca,\ca)$ we take $\mbox{\it PDiff}(M)[[\hbar]]$
($\mbox{\it PDiff}(M)$ are the polydifferential operators on
$M$),\footnote{to be precise about the tensor product $\Ga(\bw
T^*B)\otimes \mbox{\it PDiff}(M)$: the coefficients of these
differential forms on $B$ with values in polydifferential
operators on $M$ are allowed to be arbitrary smooth functions on
$M\times B$; this will also allow us to pass from $M\times B$ to
an open subset $U\ss M\times B$} $\chi$ itself is allowed to
depend formally on $\hbar$; moreover, we require that $\ga_0$ is
the ordinary product on $M$ plus $O(\hbar)$ (i.e.~a family of
$*$-products on $M$) and that $\ga_1$ is $O(\hbar)$, so that
parallel transport (and also $a_{D,x}$'s giving the inner
automorphisms) make sense. Notice that the parallel transport is
of the form $f\mapsto f+\hbar {\cal D}_1f+\hbar^2{\cal
D}_2f+\dots$ where ${\cal D}_k$'s are differential operators and
that $a_{D,x}=\exp(\int_D\ga_2)+O(\hbar)$. Just as above, for any
point $x\in B$ we have a $*$-product on $M$, parallel transport
gives isomorphisms between them and for closed contractible curves
these automorphisms are inner, given by $a_{D,x}$'s. We also have
to restrict the Lie algebra acting by  inner transformations to
those elements of total degree 1 whose (0,1)-bihomogeneous part is
$O(\hbar)$.

We shall need a slight generalization: we choose an open subset
$U\ss M\times B$ and suppose that $\ga$ is defined just on $U$
(one could also consider a submersion $E\map B$ with a transversal
foliation). The difference is that the $*$-products, parallel
transport and $a_{D,x}$'s are defined only locally now: For any
$x\in B$ let $U_x=U\cap (M\times\{x\})$; it is an open subset of
$M$. Then $\ga_0$ is a family of $*$-products on $U_x$'s. If $c$
is a curve in $B$ connecting 2 points $x_{1,2}$ and $V\ss M$ an
open subset contained in $U_x$ for every $x\in c$, then $\ga_1$
gives an isomorphism of the $*$-products on $V$ over $x_1$ and
$x_2$. Finally, if $x_1=x_2$ and $D$ is a disk in $B$ with
boundary $c$, and if moreover $V\ss U_x$ for every $x\in D$, then
we get a function $a_{D,x_1}$ on $V$ (formally depending on
$\hbar$) such that the automorphism given by $c$ is the inner
automorphism given by $a_{D,x_1}$. Moreover, if we choose a different
$D'$ homotopic to $D$ rel boundary, such that $V\ss U_x$ for every
point swept by the homotopy, then $a_{D'}=a_{D}\exp{\int\phi}$,
where the integral is over the volume swept by the homotopy.

\section{Tight families of Poisson structures and their quantization}

Let again $M$ and $B$ be manifolds and $\chi$ a closed 3-form on
$B$. Let $p_{M,B}$ be the projections from $M\times B$ to $M$ and
$B$. By a {\zv $p_B^*\chi$-tight family of Poisson structures on
$M$ indexed by $B$} we shall mean an element $\si\in\Ga(\bw
T^*B)\ot\Ga(\bw TM)$ of total degree 2, satisfying the modified
Maurer--Cartan equation (that again becomes ordinary Maurer-Cartan
equation if we augment $\Ga(\bw TM)$)
$$d\si+{1\over2}[\si,\si]=\chi,$$ where $d$ is the differential on
$B$ and $[\;,\;]$ is the Schouten bracket on $M$. To be precise
about the tensor product: $\si$ is a section of $\bw^2(p^*_B
T^*B\oplus p^*_M TM)$ where $p_{M,B}$ are the projections from
$M\times B$ to $M$ and $B$. More generally, we shall consider an
open subset $U\ss M\times B$ and suppose $\si$ is defined only
there.

Equivalently, $\si$ is a $p_B^*\chi$-Dirac structure (see
\cite{ja,sw} or the appendix) on $M\times B$ transversal to
$TM\oplus T^*B$ (this works also for arbitrary submersions $E\map
B$, not just for $M\times B\map B$; notice also that we have a
pair of transversal Dirac structures, i.e.~a Lie bialgebroid).

When $\si$ is decomposed to its bihomogeneous components
$\si=\si_0+\si_1+\si_2$ (where $\si_i$ is an $i$-form on $B$ with
values in $2-i$-vectors on $M$) then we get the same equations as
for $\ga$ above, with a similar interpretation: just replace
isomorphisms of algebras by Poisson diffeomorphisms and inner
automorphisms by Hamiltonian diffeomorphisms.

We again have an action of 2-forms on $B$ by $\si\mapsto\si+\be$,
making $\si$ to a $\chi+d\be$-tight family (the {\zv outer
transformations}), and of the Lie algebra of the elements of total
degree 1 (the {\zv inner transformations}) by
$\al\cdot\si=d\al+[\al,\si]$ (these $\al$'s are sections of the
complementary Dirac structure).

For the purpose of this note, let a {\zv formal $p_B^*\chi$-tight
family of Poisson structures} be a $\si$ as above which is a
formal power series in $\hbar$ and moreover $\si_0$ and $\si_1$
are $O(\hbar)$ (again, $\chi$ may also be allowed to depend on
$\hbar$). For the inner transformations we need a little
restriction, the (0,1)-bihomogeneous component of $\al$ must be
$O(\hbar)$.

{\em As an easy application of the Formality theorem, such a $\si$
can always be quantized to a $\chi$-tight family of $*$-products}
(see the last section). Moreover, the Formality theorem gives a
bijection between their equivalence classes with respect to the
inner transformations, and the quantization is equivariant
w.r.t.~the outer transformations.

Let now $\psi$ be a closed 3-form on $M\times B$; a {\zv
$\psi$-tight family of Poisson structures} is a $\psi$-Dirac
structure  on $M\times B$ transversal to $TM\oplus T^*B$ (again it
can be expressed as an element of $\Ga(\bw T^*B)\ot\Ga(\bw TM)$ of
total degree 2, satisfying some equation not written here). If
$\be$ is a 2-form on $M\times B$ and $\si$ a $\psi$-tight family,
we can define a $\psi+d\be$-Dirac structure $\tau_\be\si$. It need
not be transversal to $TM\oplus T^*B$; however, if $\si$ is a
formal tight family, $\tau_\be\si$ is also a formal tight family.
 If there is a
2-form $\be$ on $M\times B$ such that $\psi+d\be=p_B^*\chi$, we
can quantize $\tau_\be\si$ to a $\chi$-tight $*$-product family;
this can be called a quantization of $\si$ itself.

\section{Quantization of twisted Poisson structures to tight
$*$-product families}

Let $\phi$ be a closed 3-form on $M$ (possibly depending formally
on $\hbar$) . We choose $B=M$ (other choices would do in special
cases). Let $U\ss M\times B$ be an open subset containing the
diagonal such that $p_M^*\phi-p_B^*\phi$ is exact on $U$ (e.g. $U$
is a tubular neighborhood of the diagonal).

Let $\pi$ be a formal $\phi$-Poisson structure on $M$; we make it
to a (constant) $p_M^*\phi$-tight Poisson family (by setting
$\si_0=\pi$, $\si_1=\si_2=0$). We restrict this family from
$M\times B$ to $U$ and denote it $\si$.

Now choose a 2-form $\be$ on $U$ s.t.~$p_B^*\phi=p_M^*\phi+d\be$.
Then $\tau_\be\si$ is a $p_B^*\phi$-tight Poisson family, hence it
can be quantized to a $\phi$-tight family of $*$-products.

\section{Quantization of twisted Poisson structures to stacks}

Suppose now that the periods of $\phi$ are in $\dvapiz$ (hopefully
this $\pi=3.14\dots$ will not be confused with the bivector field
$\pi$). Let us also choose a cohomology class $\Phi\in
H^3(M,\dvapiz)$ whose image in $H^3(M,{\mathbb C})$ is the class
of $\phi$.

Let us choose a triangulation of $M$ which is fine enough (or $U$
is large enough) for the following property to hold: if $x$ is in
a simplex $s$ then $S_s\ss U_x$, where $S_s$ is the star of $s$.
Recall that the star is the union of the interiors of the
simplices containing $s$; $S_s$ is the intersection of $S_i$'s,
where $i$'s are the vertices of $s$, and $S_i$'s, where $i$'s are
all the vertices of the complex, form an open cover of $M$:
$$\epsfxsize=4cm \epsfbox{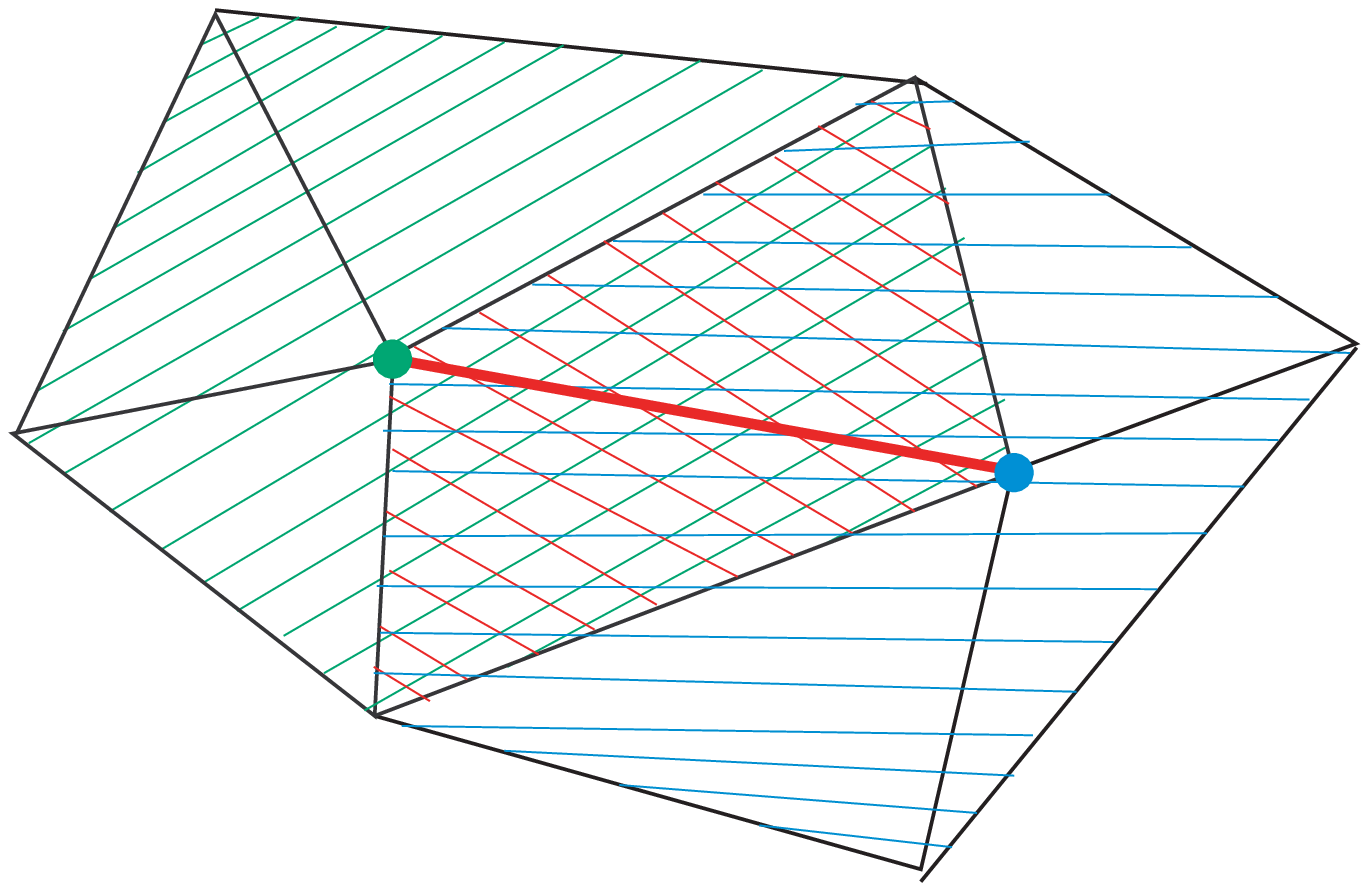}$$

Whenever two vertices $i$, $j$ form an edge, we denote $T_{ij}$
the parallel transport from $i$ to $j$ along the edge, whenever
$i$, $j$ and $k$ form a triangle $\triangle$, we denote $a_{ijk}$
the function $a_{\triangle,i}$, and whenever $i$, $j$, $k$ and $l$
form a tetrahedron, we denote $c_{ijkl}=\exp(\int\phi)$, where we
integrate over the tetrahedron. On each $S_i$ we have a $*$
product, $T_{ij}$ gives isomorphisms on the overlaps $S_i\cap S_j$
and on triple overlaps these isomorphisms are inner, given by
$a_{ijk}$'s. We obviously have the identities
$$T_{ij}=T_{ji}^{-1},\; a_{ijk}=a_{ikj}^{-1},\;
a_{jki}=T_{ij}(a_{ijk}),\; T_{ki}T_{jk}T_{ij}=\mbox{\it
Ad}(a^{-1}_{ijk})$$ and
$$a_{ikl}\,a_{ijk}=c_{ijkl}\,a_{ijl}\,T_{ji}(a_{jkl})$$ (the
products between $a$'s are the $*$-products in $S_i$).

Now we would like to multiply each $a_{ijk}$ by a non-zero number
$b_{ijk}$ so as to remove the factors $c_{ijkl}$ from the last
formula; this is indeed possible, since the periods of $\phi$ are
in $\dvapiz$, hence $c_{ijkl}$ is a coboundary. The choice of
$b$'s is fixed modulo coboundaries by the cohomology class $\Phi$.

According to Kashiwara \cite{kashi}, the $*$-products on $S_i$'s,
the isomorphisms $T_{ij}$'s and the modified $a_{ijk}$'s define a
stack of algebras on $M$.

\section{The proof involving Formality theorem}

In this section we will show how Formality theorem may be used to
quantize tight Poisson families to tight $*$-product families. We
will use the language of graded $Q$-manifolds, as in \cite{def}.
Recall that Kontsevich constructed there an equivariant map
between the graded $Q$-manifolds ${\cal U}:\Ga(\bw
TM)[2]\map\mbox{\it PDiff}(M)[2]$, sending 0 to the ordinary
product $m_0$ (the map is given by a divergent Taylor series, but
that is enough for deformation quantization; {\it PDiff} denotes
the space of polydifferential operators); moreover, the graded
$Q$-manifold $\mbox{\it PDiff}(M)[2]$ (more precisely, the formal
neighborhood of $m_0$) becomes isomorphic to the direct product of
$\Ga(\bw TM)[2]$ with a contractible graded $Q$-manifold.

For our purpose it is important that $\cal U$ is also equivariant
w.r.t.~translations by constants, and hence it can be extended to
the augmented spaces. A grading-and-$Q$-equivariant map from
$T[1]B$ to the augmented $\Ga(\bw TM)[2]$ is the same as a choice
of a closed 3-form $\chi$ on $B$ and of a $p_B^*\chi$-tight
Poisson family $\si$. We will allow this map to depend formally on
$\hbar$ in such a way that setting $\hbar=0$ yields a map constant
to second order along $B\ss T[1]B$; this is a formal tight Poisson
family. Now the composition with the augmented $\cal U$ is a
formal $\chi$-tight $*$-product family.

\appendix

\section{Exact Courant algebroids and $\phi$-Dirac structures}

A {\zv Courant algebroid} over a manifold $M$ is a vector bundle
$E\to M$ equipped with a field of nondegenerate  symmetric
bilinear forms
 $( \cdot , \cdot )$ on the fibres,  an
$\mathbb R$-bilinear bracket $[\cdot , \cdot ]:
\Gamma(E)\times\Gamma(E)\rightarrow\Gamma(E)$
 on the space of sections of $E$,
and a bundle map $\rho :E\to TM$ (the {\zv anchor}),
 such that the following
properties are  satisfied:
 \begin{enumerate}
\item for any $e_{1}, e_{2}, e_{3}\in \Gamma (E)$,
$[e_{1},[e_{2}, e_{3}]]=[[e_1,e_2],e_3]+[e_2,[e_1,e_3]];$
\item  for any $e_{1}, e_{2} \in \Gamma (E)$,
$\rho [e_{1}, e_{2}]=[\rho e_{1}, \rho  e_{2}];$
\item  for any $e_{1}, e_{2} \in \Gamma (E)$ and $f\in C^{\infty} (M)$,
$[e_{1}, fe_{2}]=f[e_{1}, e_{2}]+(\rho (e_{1})f)e_{2} ;$
\item for any $e, h_{1}, h_{2} \in \Gamma (E)$,
  $\rho (e) (h_{1}, h_{2})=([e , h_{1}] ,
h_{2})+(h_{1}, [e , h_{2}] )$;
\item for any $e,h\in\Gamma(E)$, $(e,[h,h])=([e,h],h)$.
\end{enumerate}

Equivalently, instead of the bracket $[\cdot,\cdot]$, we can use a
linear map $e\mapsto Z_e$ which maps sections of $E$ to vector
fields on the total space of $E$. The vector field $Z_e$ is a lift
of $\rho(e)$ from $M$ to $E$, and the first four axioms just say
that the flows of $Z_e$'s preserve the structure of $E$. The
bracket $[e_1,e_2]$ is the Lie derivative of $e_2$ by $Z_{e_1}$.

A {\zv Dirac structure} in $E$, also called an $E$-Dirac structure
on $M$, is a maximal isotropic subbundle $L$ of $E$ whose sections
are closed under the bracket, i.e. which is  preserved by the flow
of $Z_e$ for any $e\in \Gamma(L)$. The restriction of the bracket
and anchor to any Dirac structure $L$ form a Lie algebroid
structure on $L$.

It follows from the definition that the sequence $0\map T^*M\map
E\map TM\map 0$ is a complex (the second arrow is $\rho^*$
composed with the isomorphism $E^*\simeq E$ given by
$(\cdot,\cdot)$, the third arrow is $\rho$). If it is an exact
sequence, $E$ is an {\zv exact Courant algebroid (ECA)}. They are
classified by the 3rd de Rham cohomology \cite{ja} (see also
\cite{sw}): If we choose an isotropic splitting of the exact
sequence, so that $E$ becomes $TM\oplus T^*M$ with the bilinear
form $((X_1, \xi_1),( X_2,\xi_2))=\xi_1(X_2)+\xi_2(X_1)$, the
bracket is
$$
[(X_1, \xi_1),( X_2,\xi_2)]= ([X_1,X_2], {\cal
L}_{X_1}\xi_2-i_{X_2}d\xi_1+\phi(X_1,X_2,\cdot))
$$
for some closed 3-form $\phi$; this ECA will be denoted $(TM\oplus
T^*M)_\phi$. The space of isotropic splitting is an affine space
over the space of 2-forms, and adding a 2-form $\be$ changes
$\phi$ to $\phi+d\be$. In other words, there is an action of the
additive group of 2-forms on the bundle $TM\oplus T^*M$ given by
$\tau_\be(X,\xi)=(X,\xi+\be(X,\cdot))$; $\tau_\be$ is an
isomorphism from $(TM\oplus T^*M)_\phi$ to $(TM\oplus T^*M)_
{\phi+d\be}$. The cohomology class of $\phi$ is the {\zv
characteristic class of the ECA}.

 A $(TM\oplus T^*M)_\phi$-Dirac structure will be called
simply a {\zv $\phi$-Dirac structure}. If $\pi$ is a bivector
field on $M$, the graph of $\tilde{\pi}:T^*M\to TM$ is a
$\phi$-Dirac structure if and only if $\pi$ is a $\phi$-Poisson
structure, i.e. if $[\pi,\pi]=2\wedge^3\ti\pi(\phi)$.

Although it is not strictly needed for the paper, I include here
graded Q-manifold version of some of the notions defined above.
Let $\Upsilon$ be a principal $\rc[2]$-bundle over $T[1]M$ in the
category of graded $Q$-manifold. Such $\Upsilon$'s are classified
by the 3rd de Rham cohomology of $M$. Indeed: there is always a
grading-preserving splitting of $\Upsilon$ to $T[1]M\times\rc[2]$,
$Q=d+\phi\,\partial/{\partial t}$, where $d$ is the de Rham
differential on $M$, $t$ the coordinate on $\rc[2]$ and $\phi$ a
closed 3-form on $M$. To pass to another splitting we have to
choose a 2-form $\be$ on $M$ and after a simple computation we
find that $\phi$ gets changed to $\phi+d\be$. There is an ECA
corresponding to such a $\Upsilon$: its sections are the vector
fields on $\Upsilon$ of degree $-1$, the bracket is
$[[Q,v_1],v_2]$ and the bilinear form is $[v_1,v_2]$.

Let now $j^1\Upsilon$ be the space of 1-jets of sections of
$\Upsilon\map T[1]M$. The Dirac structures in the corresponding
ECA can be encoded as Legendrian graded $Q$-submanifolds of
$j^1\Upsilon$. Their formal neighborhoods in the
infinite-dimensional graded $Q$-manifold of all Legendrian
submanifolds of $j^1\Upsilon$ give rise to $L_\infty$-algebras;
the graded Lie algebra $\Ga(\bw TM)$ appears in the case $\pi=0$,
$\phi=0$.

\end{document}